\documentclass[11pt]{article}

\usepackage{amsmath,amssymb}
\usepackage[margin=1.15in]{geometry}
\usepackage{mathpazo}
\usepackage{microtype}
\usepackage{booktabs}
\usepackage{xcolor}
\usepackage{tikz}
\usetikzlibrary{arrows.meta,calc,decorations.pathreplacing,shapes.geometric}
\usepackage{pgfplots}
\pgfplotsset{compat=1.18}
\usepackage{tcolorbox}
\usepackage{caption}
\usepackage[colorlinks=true,linkcolor=black,urlcolor=blue!60!black,citecolor=black]{hyperref}

\newcommand{\cl}{\ensuremath{\clubsuit}}
\newcommand{\di}{\textcolor{red!80!black}{\ensuremath{\diamondsuit}}}
\newcommand{\he}{\textcolor{red!80!black}{\ensuremath{\heartsuit}}}
\newcommand{\spd}{\ensuremath{\spadesuit}}

\newcommand{\cardat}[5]{%
  \begin{scope}[shift={(#1,#2)}]
    \draw[rounded corners=2.5pt,fill=white,line width=0.6pt] (-.38,-.55) rectangle (.38,.55);
    \node[font=\footnotesize\bfseries,text=#5] at (-.17,.32) {#3};
    \node[font=\Large,text=#5] at (.06,-.14) {$#4$};
  \end{scope}}
\newcommand{\backat}[2]{%
  \begin{scope}[shift={(#1,#2)}]
    \draw[rounded corners=2.5pt,fill=blue!55!black,line width=0.6pt] (-.38,-.55) rectangle (.38,.55);
    \draw[rounded corners=1.5pt,white,line width=0.5pt] (-.3,-.47) rectangle (.3,.47);
    \draw[white,line width=0.4pt] (-.3,-.47) -- (.3,.47) (-.3,.47) -- (.3,-.47);
  \end{scope}}

\newtcolorbox{rulebox}[1]{colback=gray!6,colframe=gray!55,boxrule=0.5pt,arc=2pt,
  left=6pt,right=6pt,top=4pt,bottom=4pt,fonttitle=\bfseries,title={#1},
  coltitle=black,colbacktitle=gray!18}

\title{\LARGE Around the Clock\\[2mm]
\large New Twists on Fitch Cheney's Five-Card Trick}
\author{%
Doron Levy\thanks{%
Department of Mathematics, University of Maryland, College Park.
Email: \texttt{dlevy@umd.edu}.}
\and
Yaakov Malinovsky\thanks{%
Department of Mathematics and Statistics,
University of Maryland, Baltimore County.
Email: \texttt{yaakovm@umbc.edu}.
Research supported in part by BSF grant 2020063.}
}
\date{\today}

\providecommand{\harvarditem}[4][]{\bibitem{#4}}
\begin{document}
\maketitle

\begin{abstract}
\noindent
Fitch Cheney's five-card trick is one of the purest two-person effects
in card magic, and it has attracted mathematicians for more than
seventy five years. This article retells the effect and its classic
method, surveys how far it can be pushed, and presents several new
results that matter to performers. We give a version of Cheney's
method that needs no fixed first suit and no fixed ``lowest card'':
shifting every card one step in rank and one step in suit shifts the display in the same way.
A heckler who may swap two neighboring cards of the display defeats every method with a full deck. 
We show that a $29$-card deck is the largest for which such a swap can be corrected, and that the corresponding maximum is $14$ cards if any two cards may be swapped. We give simple methods that achieve both bounds.
Finally, we make Berlekamp's arithmetic method for the $124$-card deck independent of where the numbering starts. We also describe a pocket version with eight cards.
\end{abstract}

\section{The effect}

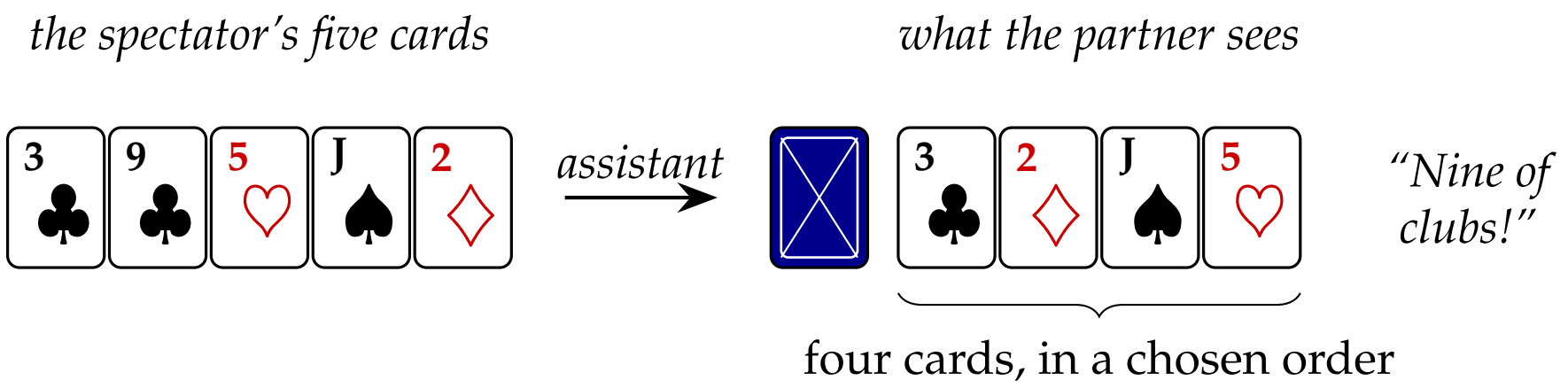
\begin{figure}[ht]
\centering
\begin{tikzpicture}[scale=0.95]
  \node[font=\small\itshape] at (1.6,1.25) {the spectator's five cards};
  \cardat{0}{0}{3}{\clubsuit}{black}
  \cardat{0.8}{0}{9}{\clubsuit}{black}
  \cardat{1.6}{0}{5}{\heartsuit}{red!80!black}
  \cardat{2.4}{0}{J}{\spadesuit}{black}
  \cardat{3.2}{0}{2}{\diamondsuit}{red!80!black}
  \draw[-{Stealth[length=3mm]},thick] (4,0) -- (5.2,0) node[midway,above,font=\small\itshape] {assistant};
  \node[font=\small\itshape] at (8.2,1.25) {what the partner sees};
  \backat{6}{0}
  \cardat{7}{0}{3}{\clubsuit}{black}
  \cardat{7.8}{0}{2}{\diamondsuit}{red!80!black}
  \cardat{8.6}{0}{J}{\spadesuit}{black}
  \cardat{9.4}{0}{5}{\heartsuit}{red!80!black}
  \draw[decorate,decoration={brace,amplitude=5pt,mirror}] (6.62,-.75) -- (9.78,-.75)
     node[midway,below=6pt,font=\small] {four cards, in a chosen order};
  \node[font=\small\itshape,align=center] at (11.1,0) {``Nine of\\clubs!''};
\end{tikzpicture}
\caption{The effect. The only signal is \emph{which} four cards are showing and
\emph{in what order}. The four cards and their order could be communicated over the telephone.}
\label{fig:effect}
\end{figure}

A spectator chooses any five cards from an ordinary deck.
She hands them to the magician's \emph{assistant}, who looks them over, turns one face down, and lays the other four face up in a chosen order.
The assistant's \emph{partner}, who has been out of the room the whole time, walks in, glances at the row, and names the face-down card (Figure~\ref{fig:effect}).

The assistant's \emph{partner}, who has been out of the room the whole time, walks in, glances at the row, and names the face-down card.

The assistant and the partner agree on a method in advance. There are no signals in what the assistant says, no glances, and no timing cues. Everything the partner learns is contained in the row itself, so the trick can be done over the phone, by text message, or with the partner in another city. This is why mathematicians love the trick, and why audiences who know a little about ``two-person codes'' find it baffling.

It is useful to separate the theatrical effect from the mathematical problem.
A \emph{hand} is an unordered set of five cards: the order in which the cards were chosen carries no information.
A \emph{row} is an ordered list of four distinct cards. A method must assign to every five-card hand an ordered row of four cards from that hand, and from the row alone the partner must be able to determine which fifth card was hidden.
In particular, two different hands must never produce the same row. This simple ``different hands need different messages'' principle is the counting idea used repeatedly below.

\section{Cheney's method}

The method rests on three simple observations.

\paragraph{Two cards share a suit.} 
With five cards and only four suits, at least two cards have the same suit. One of these two cards is hidden. The other is shown first, and it announces the suit of the hidden card.
Call the first displayed card the \emph{base}.

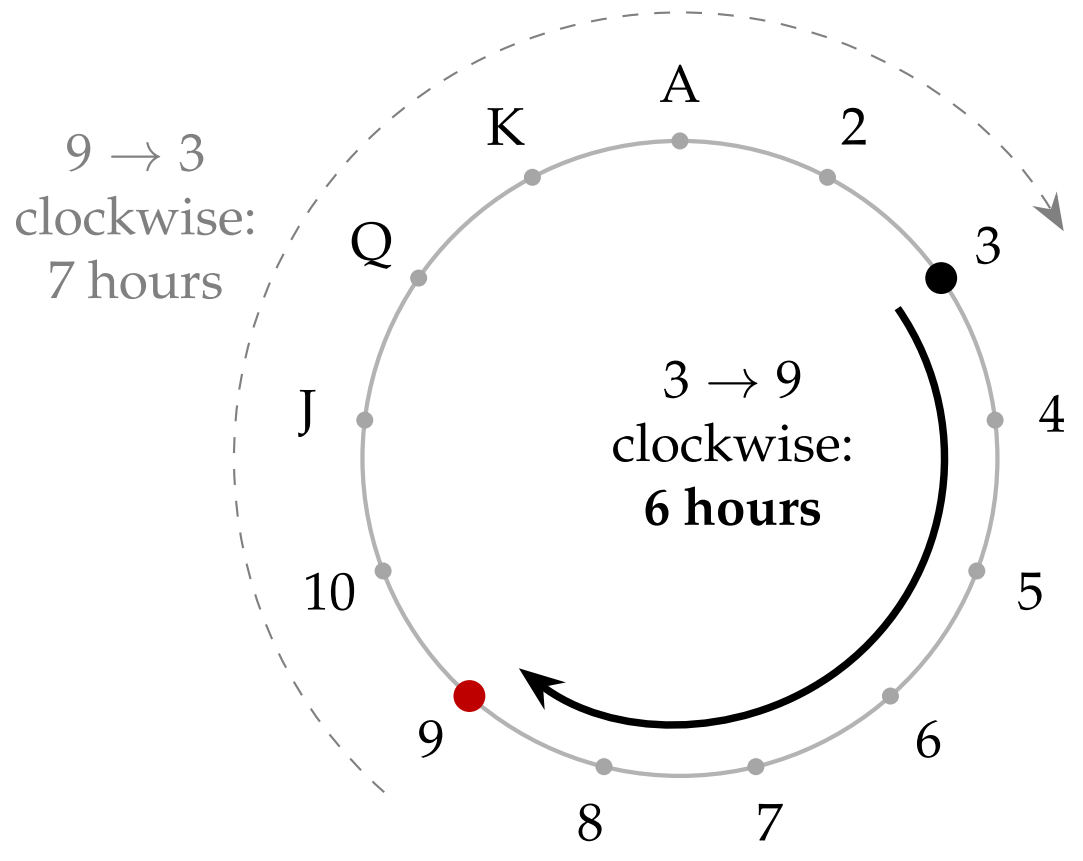
\begin{figure}[ht]
\centering
\begin{tikzpicture}[scale=1]
  \def\R{2.1}
  \draw[gray!60,line width=0.8pt] (0,0) circle (\R);
  \foreach \i/\rk in {0/A,1/2,2/3,3/4,4/5,5/6,6/7,7/8,8/9,9/10,10/J,11/Q,12/K}{
     \pgfmathsetmacro\ang{90-360*\i/13}
     \fill[gray!70] (\ang:\R) circle (1.6pt);
     \node[font=\small] at (\ang:\R+0.38) {\rk};
  }
  \pgfmathsetmacro\aa{90-360*2/13}
  \pgfmathsetmacro\ab{90-360*8/13}
  \fill[black] (\aa:\R) circle (3pt);
  \fill[red!75!black] (\ab:\R) circle (3pt);
  \draw[-{Stealth[length=3mm]},line width=1.4pt,black] (\aa:\R-0.35) arc[start angle=\aa,end angle=\ab+4,radius=\R-0.35];
  \node[font=\small,align=center] at (0.35,0.1) {$3 \to 9$\\clockwise:\\\textbf{6 hours}};
  \draw[-{Stealth[length=2.5mm]},dashed,gray] (\ab:\R+0.85) arc[start angle=\ab,end angle=\aa-360-4,radius=\R+0.85];
  \node[font=\small,gray,align=center] at (-3.6,1.6) {$9 \to 3$\\clockwise:\\7 hours};
\end{tikzpicture}
\caption{The clock of ranks. From $3$ to $9$ is six hours clockwise, from $9$ to $3$
is seven. So with $3\cl$ and $9\cl$ in the hand, $3\cl$ is the base and $9\cl$ is
hidden. The distance to signal is always between $1$ and $6$.}
\label{fig:clock}
\end{figure}

\paragraph{The clock of ranks.} Put the thirteen ranks around a clock face, Ace through King, with King followed again by Ace (Figure~\ref{fig:clock}). Start at one of the two same-suit ranks and count clockwise to the other. Suppose this takes $k$ hours, where $1\leq k\leq 12$. Going clockwise from the second rank back to the first then takes $13-k$ hours. Exactly one of $k$ and $13-k$ is at most $6$, because they add to $13$. The assistant chooses as base the card from which the other card is $1$ to $6$ hours clockwise. The other card is hidden, and this distance is the number to be communicated.

\paragraph{Three cards say a number from 1 to 6.} The partner now knows the hidden card's suit and that it lies $1$ to $6$ hours clockwise of the base. The remaining three cards can be laid out in $3!=6$ orders, exactly enough. The performers agree on a way to rank any three cards as \emph{low}, \emph{middle} and \emph{high}, and on a chart matching each order to a number. Figure~\ref{fig:code} shows the chart used in this article. Any fixed one-to-one assignment of the six orders to the numbers $1$ through $6$ works.

\begin{figure}[ht]
\centering
\begin{tikzpicture}[x=1cm,y=1cm]
  \foreach \nn/\ua/\ub/\uc [count=\k] in {1/L/M/H,2/L/H/M,3/M/L/H,4/M/H/L,5/H/L/M,6/H/M/L}{
    \pgfmathsetmacro\x{2.4*(\k-1)}
    \node[font=\bfseries] at (\x+0.5,0.85) {\nn};
    \foreach \t/\j in {\ua/0,\ub/1,\uc/2}{
      \draw[rounded corners=2pt,fill=white,line width=0.6pt] (\x+0.5*\j-0.21,-0.32) rectangle (\x+0.5*\j+0.21,0.32);
      \node[font=\small\bfseries] at (\x+0.5*\j,0) {\t};
    }
  }
\end{tikzpicture}
\caption{The code chart: the six orders of the low (L), middle (M) and high (H) cards,
in dictionary order, with $L<M<H$.}
\label{fig:code}
\end{figure}

\paragraph{An example.} The hand in Figure~\ref{fig:effect} is
$3\cl, 9\cl, 5\he, J\spd, 2\di$. The clubs are the pair, $3\cl$ is the base, $9\cl$ is hidden, and the number to send is $6$. In the classic method, cards are ranked first by rank, with Ace low, and then, in case of a tie, by suit in a fixed order. Here we use clubs, diamonds, hearts, spades, in that order. The other three cards are low $2\di$, middle $5\he$, high $J\spd$. Order $6$ is high, middle, low, so the classic display is
\[
  3\cl,\quad J\spd,\quad 5\he,\quad 2\di .
\]
The partner sees the base $3\cl$, reads the number $6$ from the order of the other three cards, counts six hours clockwise from $3$, and names the nine of clubs.

Figure~\ref{fig:effect} shows a \emph{different} row for the same hand. It comes from the new method of Section~\ref{sec:rotproof}, and the difference between the two methods is the heart of this article.

\section{A short history}

The trick is due to William Fitch Cheney Jr., a mathematician and magician,
and appeared under the title ``Telephone Stud'' in Wallace Lee's booklet
\emph{Math Miracles} in 1950 \cite{Lee}.
Michael Kleber's 2002 article ``The Best Card Trick'' in \emph{The Mathematical Intelligencer} \cite{Kleber} explains the trick and asks how large the deck can be. It also gives an arithmetic method for a 124-card deck, which Kleber learned from Elwyn Berlekamp.
The correspondence collected in \cite{MathPuzzle}
records that Jim Yingst posted essentially the same construction to the
newsgroup \texttt{rec.puzzles} in 1997.
Simonson and Holm give a detailed treatment of the trick
in the context of teaching discrete mathematics \cite{SimonsonHolm}. Colm Mulcahy's
\emph{Mathematical Card Magic} \cite{Mulcahy} develops many variants,
and a 2024 paper by a group in MIT's PRIMES STEP program studies further
variants and introduces new ones \cite{PRIMES}. The new results
presented here are developed in a technical paper written in parallel
with this article \cite{Symmetry}.

The trick also appears in many books. Martin Gardner mentions it in
passing in \emph{Mathematics, Magic and Mystery}
\cite[p.~32]{Gardner}.
Karl Fulves gives a related telephone trick, ``The Ultra Code,'' as No.~31 in \emph{More Self-Working Card Tricks} \cite[pp.~41--43]{Fulves}. It uses five cards of one suit and an extra verbal code.
Van Lint and Wilson present the trick in Example~5.1 of \emph{A Course in Combinatorics} \cite[pp.~44--45]{vanLintWilson}, in Chapter~5, ``Systems of distinct representatives,'' and prove that a 124-card deck is possible.
Mulcahy's article ``Fitch Cheney's Five Card Trick'' is reprinted in
\emph{The Edge of the Universe} \cite[pp.~273--276]{MulcahyEdge}.
Chapter~13 of his book \cite{Mulcahy}, ``The Hidden Value of Cards,''
opens with ``Fitch Cheney's Five-Card Twist'' and ends with the
four-card version ``Fitch Four Glory,'' which allows some displayed
cards to be face down.

The MIT text \emph{Mathematics for Computer Science} treats the trick
in Section~11.11, ``A Magic Trick'' \cite[pp.~346--353]{MCS}.
Devadas's \emph{Programming for the Puzzled} uses it as Puzzle~3,
``You Can Read Minds (With a Little Calibration)'' \cite{Devadas}.
Winkler's \emph{Mathematical Puzzles} includes it as ``Missing Card,''
with its solution in Chapter~13, ``Information, Please''
\cite[pp.~167--169]{Winkler}. Cl\'ement and Guillaume Deslandes pose
the 124-card version directly, with the cards numbered 0 to 123,
as Riddle~80, ``A Magic Trick,'' in \emph{The Price of Cake}
\cite{Deslandes}. The trick also appears in the recent book \emph{Mathematical Puzzles and Curiosities} by David, Khovanova and Shpilman \cite{DKS26}.

\section{How big can the deck get?}
\label{sec:ceiling}

A full deck is surprisingly roomy for this trick. The assistant can lay down $52 \times 51 \times 50 \times 49 = 6{,}497{,}400$ different rows, while the spectator can deal only $\binom{52}{5}=2{,}598{,}960$ different hands. How far can the deck grow before the rows run out? This counting ignores the restriction that the four displayed cards must belong to the hand. Therefore it can prove impossibility when there are too few rows, but by itself it cannot prove that a method exists when there are enough.

With a deck of $d\geq5$ cards, choosing the four displayed cards and then ordering them gives $4!\binom{d}{4}=24\binom{d}{4}$ possible rows, while there are $\binom{d}{5}=\binom{d}{4}\frac{d-4}{5}$ possible hands. The ratio of rows to hands is therefore
\[
  \frac{24\binom{d}{4}}{\binom{d}{5}}=\frac{120}{d-4}.
\]
It is $2.5$ for a full deck and falls to exactly $1$ when $d=124$. With $125$ or more cards there are more hands than rows, so two different hands must produce the same row, and the trick is impossible. With $124$ cards it is possible, but only just: there are exactly $225{,}150{,}024$ rows and the same number of hands, so every row must be used exactly once. Kleber gives a proof that such a perfect pairing exists \cite{Kleber}. The nontrivial point is that this can be done while assigning each hand a row made from four of its own cards.

The same count applies to a hand of $n$ cards with $n-1$ of them shown. The ratio of rows to hands is $n!/(d-n+1)$, so the trick is possible only when $d\leq n!+n-1$. This gives $8$ cards for $n=3$, $27$ for $n=4$, and $124$ for $n=5$. Section~\ref{sec:variations} shows that the two smaller bounds are also attained.

\paragraph{Berlekamp's recipe.} For a $124$-card deck, with cards numbered $0$ to $123$, Berlekamp's method can be illustrated in four steps. For the hand $\{7,23,51,88,119\}$:
\begin{enumerate}
\item Add the five numbers: $288$. The remainder on division by $5$ is $3$. Numbering the positions from $0$ in increasing order, hide the card in position $3$, namely $88$.
\item Add the other four numbers: $200$. Let $r$ be the smallest nonnegative number that must be added to reach a multiple of $5$. Here $r=0$.
\item 
List the $120$ cards that are not displayed, the \emph{unseen} cards, in increasing order, with positions numbered from $0$. Three displayed cards lie below $88$, so $88$ is in position $88-3=85$.
\item Compute $(85-r)/5=17$ and lay the four cards in order number $17$ of the $24$ possible orders, counted in dictionary order from $0$. The row is $51,\ 119,\ 23,\ 7$.
\end{enumerate}
The partner finds that the displayed order has number $17$. The four cards sum to $200$, so $r=0$. The hidden card is therefore in position $5\cdot17+0=85$ among the unseen cards, which is the $88$. The method works, but it requires some mental arithmetic.

Here is why the arithmetic works. Write the five cards in increasing order as $x_0<x_1<x_2<x_3<x_4$, and let $T$ be their total. If $T$ leaves remainder $j$ upon division by $5$, Berlekamp's remainder rule hides $x_j$. After the other four cards are displayed, let their sum be $S=T-x_j$, and let $r\in\{0,1,2,3,4\}$ be the amount that must be added to $S$ to reach a multiple of $5$. 
There are $x_j$ cards below $x_j$ in the deck, and exactly $j$ of them are displayed, so the hidden card is in position $y=x_j-j$ among the unseen cards. Here and below, $a\equiv b\pmod m$ means that $a$ and $b$ leave the same remainder when divided by $m$.
Since $j$ and $T$ have the same remainder,
\[
  y=x_j-j\equiv x_j-T=-(T-x_j)=-S\pmod5.
\]
Thus $y$ has the same remainder as $r$, and since $0\leq y\leq119$, we can write $y=r+5q$ for exactly one $q\in\{0,1,\ldots,23\}$. The $24$ values of $q$ are encoded by the $4!=24$ orders of the displayed cards. 
The partner reverses this calculation: she computes $S$ and $r$, reads $q$ from the order, and names the card in position $r+5q$ among the $120$ unseen cards.

\section{Variations}
\label{sec:variations}

The two versions below are specializations of Berlekamp's
remainder construction described above. Chen et al.\
\cite[Sections~2.3.1--2.3.2, Theorem~3 and Table~3]{PRIMES}
discuss the general construction and list the optimal deck
sizes $8$ and $27$ for hands of three and four cards,
respectively. We give explicit rules and examples for these
small cases to illustrate the method.

\subsection*{A pocket version with eight cards}

Take eight cards valued $1$ through $8$, for example the Ace to eight of one suit. The spectator picks any three and hands them to the assistant, who hides one and shows two in a chosen order.
There are exactly $\binom{8}{3}=56$ possible hands and $8 \times 7 = 56$ possible rows, so, by the count above, eight cards is the largest possible deck.
Berlekamp's recipe, specialized to three cards, gives a method.

\begin{rulebox}{The assistant}
Add the three values. If the total leaves remainder $0$ on division by $3$, hide the lowest card. If it leaves remainder $1$, hide the middle card. If it leaves remainder $2$, hide the highest card. Now look at the six cards not being shown, including the hidden card, in increasing order, and number their positions $1,2,\ldots,6$. If the hidden card is among the first three positions, show the other two cards in increasing order. If it is among the last three positions, show them in decreasing order.
\end{rulebox}

\begin{rulebox}{The partner}
Add the two values showing. List the six values not showing in increasing order, and number their positions $1,2,\ldots,6$. Find the two positions $p$ for which $p+{}$(sum of the two displayed cards) is divisible by $3$. They are one of the pairs $\{1,4\}$, $\{2,5\}$, $\{3,6\}$. If the two cards showing are in increasing order, name the card in the smaller of the two positions. If they are in decreasing order, name the card in the larger position.
\end{rulebox}

For example, consider the hand $2,5,7$. Since $2+5+7=14\equiv2\pmod3$, the assistant hides the highest card, $7$. The six unseen cards are $1,3,4,6,7,8$, and $7$ occupies position $5$, so the assistant shows $5,2$. The partner sees $5,2$, whose sum is $7$. She seeks positions with $7+p\equiv0\pmod3$, namely $2$ and $5$, containing the cards $3$ and $7$. Since the row is decreasing, she names $7$.

Why it works: write the hand as $x_0<x_1<x_2$, with total $T$. If $T$ leaves remainder $j$ on division by $3$, then $x_j$ is hidden. Let $S$ be the sum of the two displayed cards. There are $x_j-1$ values below $x_j$, and exactly $j$ of them are displayed. So among the six unseen cards, numbered from $1$, the hidden card is in position $p=x_j-j$.
Hence $p\equiv x_j-T=-S\pmod3$, so $p+S$ is divisible by $3$. Each remainder occurs exactly twice among the six positions, once in each half, so the displayed sum leaves two possible hidden cards, and the increasing or decreasing order of the row picks one.

\subsection*{Four cards from twenty-seven}

A four-card hand leaves three cards to display. For $d\geq4$ there are
$3!\binom d3$ rows and $\binom d4$ hands, with ratio $24/(d-3)$, which
is $1$ at $d=27$. Number the cards $0,1,\ldots,26$. There are
$\binom{27}{4}=27\times26\times25=17{,}550$ hands and rows. The same
steps work with $4$ in place of $5$. Add the four values. If the total
leaves remainder $j\in\{0,1,2,3\}$ on division by $4$, hide the card
in position $j$ of the sorted hand, counting from $0$. Let $S$ be the
sum of the three displayed cards and let $r\in\{0,1,2,3\}$ satisfy
$S+r\equiv0\pmod4$. Exactly as before, the position $y$ of the hidden
card among the $24$ unseen cards satisfies $y+S\equiv0\pmod4$, so it
is one of the six positions $r, r+4, \ldots, r+20$. 
The $3!=6$ orders
of the displayed cards, read with the chart of Figure~\ref{fig:code}
and ranked by value, select these six positions in increasing order. For $q=0,\ldots,5$, order number $q+1$ in
Figure~\ref{fig:code} signals unseen position $r+4q$.

\section{The heckler}
\label{sec:heckler}

Suppose that before the partner sees the row, a heckler may swap two neighboring cards in it, or leave it alone, and the partner is not told which. Can the trick survive?

Write $1,2,3,4$ for the four displayed cards from lowest to highest in the agreed ordering, so that the row on the table is one of the $4!=24$ orders of 1234. Two orders are \emph{neighbors} if one is obtained from the other by swapping two neighboring cards, and every order has exactly three neighbors. The \emph{distance} between two orders is the least number of such swaps that turns one into the other (the Kendall~$\tau$ distance~\cite{Kendall38}). For example, 1234 and 1324 are at distance~1, 1234 and 2143 at distance~2, and 1234 and 2341 at distance~3, since the card $1$ must move past the other three.

For each set of four displayed cards, the performers agree in advance
on the orders the assistant may use. We call these orders \emph{usable}.
The
partner sees either the usable order the assistant laid down or one of
its three neighbors. Think of each usable order as \emph{owning} itself and its three neighbors.
The partner decodes the row she sees by finding its owner. This works exactly when any two usable orders are at distance at least~3:
\begin{enumerate}
\item At distance~1, the partner cannot distinguish one order left unchanged from the other after a swap.
\item At distance~2, some row lies one swap from each. For example, 2134 is one swap from 1234 and one swap from 2143.
\item If all distances are at least~3, no row is within one swap of two usable orders, since those orders would then be within two swaps of each other.
\end{enumerate}
We call a set of orders whose pairwise distances are all at least~3 \emph{heckler-safe}.

Each order of a heckler-safe set owns four rows, and no row
has two owners, so the set has at most $24/4=6$ orders.
We now show directly that six orders are impossible.

Suppose there were six usable orders. Their four-row sets
would partition all $24$ orders. For $j=1,2,3,4$, let $c_j$
be the number of usable orders in which card $1$ occupies
position $j$. Among all $24$ orders, exactly $3!=6$ have
card $1$ in any specified position.
For example, a usable order with card $1$ first owns three rows in which that card is still first: the order itself and the two neighbors obtained by swapping the cards in positions $2$ and $3$, or in positions $3$ and $4$.
A usable order with card $1$ second owns exactly one row
with that card first. Orders with card $1$ third or fourth
own no such row.
This gives the first equation:
\[
 3c_1+c_2 = 6.
\]
Counting the rows with card $1$ in positions $2$, $3$ and $4$ in the same way gives
\[
\begin{aligned}
c_1+2c_2+c_3 &= 6,\\
c_2+2c_3+c_4 &= 6,\\
c_3+3c_4 &= 6.
\end{aligned}
\]
The equations imply
\[
c_1=c_2=c_3=c_4=\frac32,
\]
which is impossible because the $c_j$ are integers.
Thus a heckler-safe set has at most five orders, and
Figure~\ref{fig:heckler} shows one with five.

In coding-theory language, the hypothetical six-order set
would be a \emph{perfect single-error-correcting code}.
For four objects, the argument above is a direct proof of the nonexistence result of Buzaglo and Etzion \cite[Theorem~3]{BE15}.
Related permutation codes have been studied for flash
memories that store data as rankings of electric
charges~\cite{JMSB09}, where a small drift in charge can
swap two neighboring ranks~\cite{JSB10}.

Now count as in Section~\ref{sec:ceiling}. Suppose that at most $m$ orders are usable for each set of four cards. The heckler never changes which four cards are shown, so a deck of $d$ cards offers at most $m\binom{d}{4}$ messages for its $\binom{d}{5}=\binom d4\frac{d-4}5$ hands. This forces $d\leq5m+4$, and $m\leq5$ gives $d\leq29$. \emph{In particular, the challenge is impossible with a full deck.} We now show that the bound is attained.

\begin{figure}[ht]
\centering
\newcommand{\hkbox}[1]{{\setlength{\fboxsep}{1.2pt}\framebox[1.1em]{\strut #1}}}%
\newcommand{\hkorder}[4]{\hkbox{#1}\,\hkbox{#2}\,\hkbox{#3}\,\hkbox{#4}}%
\setlength{\tabcolsep}{0.7em}%
\begin{tabular}{ccccc}
\small 0 & \small 1 & \small 2 & \small 3 & \small 4\\[2pt]
\hkorder{1}{2}{3}{4} & \hkorder{1}{4}{3}{2} &
\hkorder{2}{3}{4}{1} & \hkorder{3}{4}{1}{2} &
\hkorder{4}{2}{1}{3}
\end{tabular}
\caption{A heckler-safe set of five orders, numbered 0 to 4.
The digits 1--4 stand for the four displayed cards from lowest to
highest. Any two of these orders are at least three neighboring
swaps apart.}
\label{fig:heckler}
\end{figure}

\paragraph{A heckler-proof trick with 29 cards.} Use Berlekamp's
recipe with the five orders of Figure~\ref{fig:heckler} in place of
all $24$. Number the cards $0$ to $28$. Sort the hand and count
positions from $0$. The assistant hides the card in position $j$,
where $j$ is the remainder of the total on division by $5$. Let $S$ be
the displayed sum and $r\in\{0,\ldots,4\}$ the amount needed to bring
$S$ to a multiple of $5$. Exactly as in Section~\ref{sec:ceiling}, the
hidden card's position among the $25$ unseen cards is $y=x_j-j\equiv
r\pmod 5$, so $y=r+5q$ for a unique $q\in\{0,\ldots,4\}$. The
assistant lays out the cards in the order numbered $q$ in
Figure~\ref{fig:heckler}. A swap does not change which four cards are
showing, so the partner computes $r$ as usual.
She then looks at the order she sees and its three neighbors. Because the five orders are heckler-safe, exactly one of these four orders appears in Figure~\ref{fig:heckler}. Its number is $q$, and she names the unseen card in position $r+5q$.

\paragraph{Example.} Take the hand $\{3,8,14,20,27\}$. The total is
$72$, with remainder $2$, so the assistant hides the $14$. The other
four add up to $58$, so $r=2$. Two displayed cards lie below $14$, so
$y=14-2=12=2+5\cdot2$ and $q=2$. Order 2 is 2341, so the row is
$8,20,27,3$. The heckler swaps the middle two cards, and the partner
sees $8,27,20,3$, that is, 2431. Of 2431 and its neighbors 4231, 2341
and 2413, only 2341 appears in Figure~\ref{fig:heckler}, so $q=2$ and
the hidden card is in position $12$ among the unseen cards
$0,1,2,4,5,6,7,9,10,11,12,13,14,\dots$, which is the $14$. 

\paragraph{Arbitrary swaps.} Suppose now that the heckler may swap any
two cards of the row, or leave it unchanged.
We show that, for five-card hands, such a heckler limits the size of the deck to $14$ cards.
As before, any two usable orders must be at least three swaps apart. If they are one swap apart, the partner cannot tell them apart. If they are two swaps apart, some row is one swap from each. But any order of four cards can be turned into any other with at most three swaps. At most one swap puts the right card in position $1$, at most one more does the same for position $2$, and at most one more for position $3$. The fourth card is then also in place.
So two usable orders must be exactly three
swaps apart. 

Recall that an order is \emph{even} or \emph{odd} according to whether
it can be obtained from 1234 by an even or odd number of swaps. This
is well defined, and each swap changes the parity.
Since every order is even or odd, among any three orders at least two have the same parity. The distance between these two
is therefore even and
at most three, hence at most two, so they cannot both be usable. Thus
at most two orders are usable. Two do suffice: 1234 and 2341 differ in
all four positions, so they are not one swap apart, and the sequence
$1234\to2134\to2314\to2341$ shows that their distance is odd, hence
three. Their possible received orders are 
$$
\begin{array}{c|l}
\text{usable order} & \text{possible received orders}\\
\hline
1234 & 1234,\ 2134,\ 3214,\ 4231,\ 1324,\ 1432,\ 1243\\
2341 & 2341,\ 3241,\ 4321,\ 1342,\ 2431,\ 2143,\ 2314
\end{array}
$$
and the two lists are disjoint. Counting hands and messages gives $\binom d5\leq2\binom d4$, which forces $d\leq14$.

The same recipe attains this bound. Number the cards $0,\ldots,13$, hide $x_j$ where $j$ is the remainder of the total mod $5$, and compute $r$ from the displayed sum. Among the ten unseen cards the hidden card has position $y=x_j-j\equiv r\pmod5$, so it is in position $r$ or $r+5$. The assistant uses order 1234 in the first case and 2341 in the second. The partner decides from the table which usable order produced the row, and names the unseen card in position $r$ or $r+5$ accordingly.

\section{A method with no starting point}
\label{sec:nostart}

Here is the idea behind the new results. The clock of ranks has no fixed starting point: King is followed by Ace. The four suits can also be placed around a small clock, clubs, diamonds, hearts, spades, and back to clubs. Now consider the move that \emph{shifts the whole deck}: every card goes up one hour in rank and one step in suit. The $3\cl$ becomes the $4\di$, the $K\spd$ becomes the $A\cl$, and so on. 
When applied to a displayed row, this shift changes every card but does not change its position in the row.

After $k$ shifts, a card's rank is back where it started exactly when $k$ is a multiple of $13$, and its suit exactly when $k$ is a multiple of $4$. Since $13$ and $4$ are relatively prime, both are back for the first time after $52$ shifts. So repeating the move from any card visits all $52$ cards, and there are $52$ different shifts, counting the one that does nothing. For example, $40$ repetitions move the ranks alone one hour, since $40\equiv1\pmod{13}$ and $40\equiv0\pmod4$, and $13$ repetitions move the suits alone one step.
Magicians will recognize a similar idea in cyclic deck arrangements such as the Si Stebbins stack.

Call a method \emph{rotation-proof} if shifting every card in the hand shifts every card in the display in the same way, without changing their positions in the row. If the hand $3\cl,9\cl,5\he,J\spd,2\di$ gives some row, then the shifted hand $4\di,10\di,6\spd,Q\cl,3\he$ must give the shifted row. In symbols, if $R$ is the one-step shift and $E(H)$ the row produced from a hand $H$, then $$E(R(H))=R(E(H)),$$ with $R$ applied to each displayed card separately.

Why should a performer care? Because a rotation-proof method treats every card of the deck alike. In the method below, we never ask which rank is lowest or which suit comes first. We only ask how cards sit \emph{relative to each other} around the two clocks. In practice this means:
\begin{enumerate}
\renewcommand{\labelenumi}{(\roman{enumi})}
\item \textbf{No first suit to memorize.} Only the agreed cyclic order of the suits matters. The method works with any fixed cyclic order, including the one in your favorite stack, provided that both performers use the same order and direction.
\item \textbf{No seam between King and Ace.} The Ace is not ``low'' and the King is not ``high.'' A hand with a King and a two is handled just like the hand obtained by moving every rank six hours forward, in which the King becomes a six and the two an eight. The row moves in the same way.
\item \textbf{Practice one hand, learn fifty-two.} Every hand comes in a family of $52$ shifted versions, all handled the same way. The $2{,}598{,}960$ hands fall into $2{,}598{,}960/52=49{,}980$ families.
\end{enumerate}

These $52$ shifted hands are all distinct. Suppose a shift other than the identity carried a five-card hand back to itself. Repeatedly applying that shift partitions the deck into cycles of a common length $m>1$ dividing $52$. For example, a one-step shift of suits alone gives cycles of length $4$. 
Since the shift carries the hand onto itself, a hand that contains a card must contain that card's whole cycle. So the hand is made of complete cycles, and $m$ divides $5$.
No integer greater than $1$ divides both $5$ and $52$. Hence only the identity shift fixes a five-card hand, and every family has $52$ members.

Cheney's classic method is \emph{not} rotation-proof, as the example
in Section~\ref{sec:rotproof} will show. The absolute ranking of cards
into low, middle and high need not respect shifts, and a fixed rule
for choosing among several same-suit pairs need not respect them
either.

\section{Cheney's method, made rotation-proof}
\label{sec:rotproof}

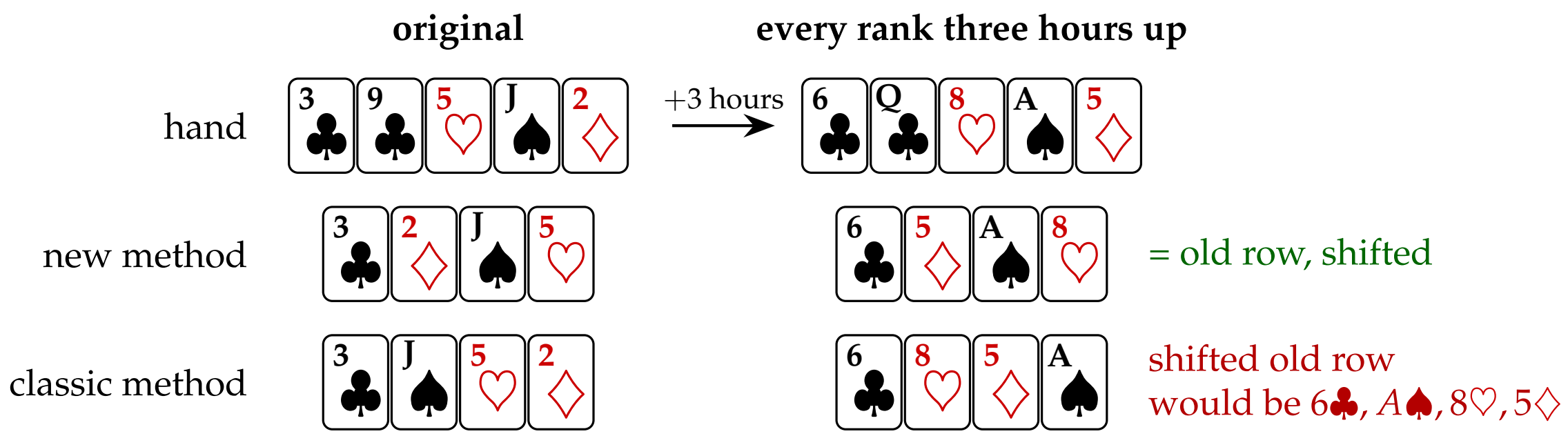
\begin{figure}[ht]
\centering
\begin{tikzpicture}[scale=0.82]
  \node[anchor=east,font=\small] at (-0.7,0) {hand};
  \node[anchor=east,font=\small] at (-0.7,-1.5) {new method};
  \node[anchor=east,font=\small] at (-0.7,-3.0) {classic method};

  \node[font=\small\bfseries] at (1.6,1.1) {original};

  \cardat{0}{0}{3}{\clubsuit}{black}
  \cardat{0.8}{0}{9}{\clubsuit}{black}
  \cardat{1.6}{0}{5}{\heartsuit}{red!80!black}
  \cardat{2.4}{0}{J}{\spadesuit}{black}
  \cardat{3.2}{0}{2}{\diamondsuit}{red!80!black}

  \cardat{0.4}{-1.5}{3}{\clubsuit}{black}
  \cardat{1.2}{-1.5}{2}{\diamondsuit}{red!80!black}
  \cardat{2.0}{-1.5}{J}{\spadesuit}{black}
  \cardat{2.8}{-1.5}{5}{\heartsuit}{red!80!black}

  \cardat{0.4}{-3.0}{3}{\clubsuit}{black}
  \cardat{1.2}{-3.0}{J}{\spadesuit}{black}
  \cardat{2.0}{-3.0}{5}{\heartsuit}{red!80!black}
  \cardat{2.8}{-3.0}{2}{\diamondsuit}{red!80!black}

  \draw[-{Stealth[length=3mm]},thick]
    (4.1,0) -- (5.3,0)
    node[midway,above,font=\scriptsize] {$+3$ hours};

  \node[font=\small\bfseries] at (7.6,1.1)
    {every rank three hours up};

  \cardat{6}{0}{6}{\clubsuit}{black}
  \cardat{6.8}{0}{Q}{\clubsuit}{black}
  \cardat{7.6}{0}{8}{\heartsuit}{red!80!black}
  \cardat{8.4}{0}{A}{\spadesuit}{black}
  \cardat{9.2}{0}{5}{\diamondsuit}{red!80!black}

  \cardat{6.4}{-1.5}{6}{\clubsuit}{black}
  \cardat{7.2}{-1.5}{5}{\diamondsuit}{red!80!black}
  \cardat{8.0}{-1.5}{A}{\spadesuit}{black}
  \cardat{8.8}{-1.5}{8}{\heartsuit}{red!80!black}

  \cardat{6.4}{-3.0}{6}{\clubsuit}{black}
  \cardat{7.2}{-3.0}{8}{\heartsuit}{red!80!black}
  \cardat{8.0}{-3.0}{5}{\diamondsuit}{red!80!black}
  \cardat{8.8}{-3.0}{A}{\spadesuit}{black}

  \node[font=\small,green!40!black,anchor=west]
    at (9.5,-1.5) {= old row, shifted};
  \node[font=\small,red!70!black,anchor=west,align=left]
    at (9.5,-3.0)
    {shifted old row\\would be $6\cl,A\spd,8\he,5\di$};
\end{tikzpicture}
\caption{Shift every rank of a hand up three hours, leaving suits
unchanged. The rotation-proof method shifts its row along with the
hand. The classic method, still correct, switches to a different
arrangement, because the $J\spd$ has wrapped past the King to become
an Ace.}
\label{fig:shift}
\end{figure}

In this section we present a rotation-proof version of Cheney's
method. Figure~\ref{fig:shift} shows what rotation-proof means in practice. Shift every rank in the first example three hours forward, leaving suits unchanged. The new method's row is the old row, shifted. 
The classic method still gives a correct row, but not the shifted copy of its old row. The reason is that the shift carries the $J\spd$ past the King, where it becomes an Ace, the lowest card.

The repair changes only two things: how the three \emph{signaling cards} (the displayed cards other than the base) are ranked, and how the same-suit pair is chosen when there is a choice.

\begin{rulebox}{Rule 1. Rank from the base.}
Rank the three signaling cards by their clockwise rank distance from the base, with fewer hours ranked lower. A card with the same rank as the base has distance zero. If two signaling cards share a rank, the one whose suit is fewer clockwise steps from the base's suit ranks lower. Use the agreed suit cycle: clubs, diamonds, hearts, spades, and back to clubs.
\end{rulebox}

To carry out Rule~1, give each card a \emph{label} relative to the base: the pair (rank distance, suit distance), where the rank distance $0,1,\ldots,12$ counts hours clockwise from the base rank and the suit distance $0,1,2,3$ counts steps clockwise from the base suit. Labels are compared by rank distance first, and by suit distance when the rank distances are equal. Different cards have different labels. A shift of the whole deck changes no label, because it moves the base and the other card by the same amount.

\begin{rulebox}{Rule 2. Choose the pair by shape.}
For each same-suit pair, find its base and hidden card as in Cheney's original method. Choose a pair with the smallest distance from base to hidden card. If several pairs tie, take each tied base in turn, label the other four cards relative to it as in Rule~1 (including the card that would be hidden), and sort the four labels from smallest to largest. Choose the tied base whose sorted list is smallest: at the first position where two lists differ, the smaller label wins.
\end{rulebox}

Once the smallest base-to-hidden distance is fixed, each tied base determines its hidden card uniquely. The partner's side of the trick is unchanged except for Rule~1: she ranks the other three cards relative to the base, reads the number from the chart, and counts that many hours clockwise from the base, keeping its suit. She does not need Rule~2.

\begin{figure}[ht]
\centering
\begin{tikzpicture}[scale=1]
  \def\R{2.2}
  \draw[gray!60,line width=0.8pt] (0,0) circle (\R);
  \foreach \i in {0,...,12}{
     \pgfmathsetmacro\ang{90-360*\i/13}
     \fill[gray!60] (\ang:\R) circle (1.4pt);
     \node[font=\scriptsize,gray] at (\ang:\R-0.3) {\i};
  }
  \node[font=\scriptsize,gray,align=center] at (0,-0.1)
    {hours\\clockwise\\from the base};

  \cardat{0}{\R+0.85}{3}{\clubsuit}{black}
  \node[font=\scriptsize\bfseries] at (0.75,\R+0.85) {base};

  \pgfmathsetmacro\ang{90-360*2/13}
  \cardat{{(\R+0.85)*cos(\ang)}}{{(\R+0.85)*sin(\ang)}}
    {5}{\heartsuit}{red!80!black}

  \pgfmathsetmacro\ang{90-360*6/13}
  \backat{{(\R+0.85)*cos(\ang)}}{{(\R+0.85)*sin(\ang)}}

  \pgfmathsetmacro\ang{90-360*8/13}
  \cardat{{(\R+0.85)*cos(\ang)}}{{(\R+0.85)*sin(\ang)}}
    {J}{\spadesuit}{black}

  \pgfmathsetmacro\ang{90-360*12/13}
  \cardat{{(\R+0.85)*cos(\ang)-0.55}}{{(\R+0.85)*sin(\ang)}}
    {2}{\diamondsuit}{red!80!black}

  \foreach \hh/\lab/\dx in {2/low/0,8/middle/0,12/high/-0.55}{
     \pgfmathsetmacro\ang{90-360*\hh/13}
     \node[font=\small\bfseries,blue!60!black]
       at ({(\R+0.85)*cos(\ang)+\dx},
           {(\R+0.85)*sin(\ang)-0.85}) {\lab};
  }

  \pgfmathsetmacro\ang{90-360*6/13}
  \node[font=\small,align=center]
    at ({(\R+0.85)*cos(\ang)},
        {(\R+0.85)*sin(\ang)-1.0}) {hidden $9\cl$};
\end{tikzpicture}
\caption{The hand $3\cl,9\cl,5\he,J\spd,2\di$ seen from its base.
The $2\di$ sits twelve hours clockwise of the $3\cl$, so it is the
\emph{high} signaling card, even though its rank is the lowest in
the usual Ace-low ordering.}
\label{fig:frombase}
\end{figure}
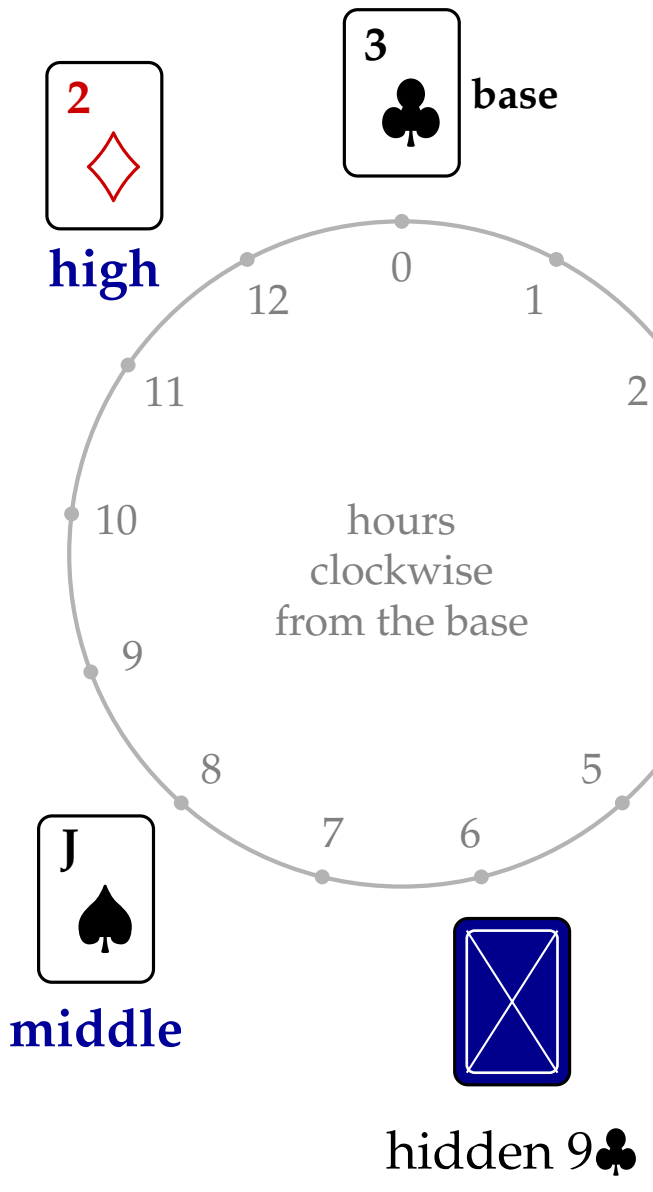

\paragraph{The first example again.}
For $3\cl,9\cl,5\he,J\spd,2\di$ the base is $3\cl$ and the number is $6$, as before. Seen from the base (Figure~\ref{fig:frombase}), the $5\he$ is two hours on, the $J\spd$ eight hours on, and the $2\di$ twelve hours on. So low, middle, high are $5\he,J\spd,2\di$, and order $6$, high-middle-low, gives the row $3\cl,\ 2\di,\ J\spd,\ 5\he$ of Figure~\ref{fig:effect}. The classic method treats the $2\di$ as the lowest signaling card. Rule~1 treats it as the highest.

\paragraph{A hand with a tie.}
The hand $A\cl,A\di,A\he,3\cl,3\he$ contains two same-suit pairs that are both two hours apart, based at $A\cl$ and at $A\he$. From $A\cl$, the nearest other card is $A\di$, with label $(0,1)$. From $A\he$, the nearest is $A\cl$, with label $(0,2)$, since the cycle runs hearts, spades, clubs. Since $(0,1)<(0,2)$, the base $A\cl$ wins, $3\cl$ is hidden, and the number is $2$. From $A\cl$, the other three rank low $A\di$, middle $A\he$, high $3\he$, and order $2$ gives $A\cl,\ A\di,\ 3\he,\ A\he$.

\paragraph{Why it always works.}
The partner does not need to know how the pair was chosen. The chosen pair consists of two cards of the same suit, oriented so that the hidden card is between one and six hours clockwise from its base, and Rule~1 and the chart communicate that distance exactly as in Cheney's original trick. Correct decoding also shows that two different hands cannot produce the same row.

It remains to show that Rule~2 selects a unique pair. Suppose two distinct candidate bases remained tied after all four labels had been compared. They would see exactly the same relative pattern of the other four cards. 
The shift that takes the first base to the second would then take each of the other four cards to the card with the same label relative to the second base. So it would carry the whole hand onto itself, which Section~\ref{sec:nostart} shows is impossible.

Finally, the quantities Rule~2 compares are all relative, so a global shift changes none of them. A shifted hand therefore chooses the shifted copy of the same pair, Rule~1 gives the shifted signaling cards the same low-middle-high pattern, and the displayed row is the old row shifted card by card.

\paragraph{What it costs.}
In most hands, little tie-breaking is needed. An enumeration of all $2{,}598{,}960$ hands shows that the pair is settled by the fewest-hours rule alone in about $85\%$ of hands. Comparing the first relative label for each tied base settles another $8\%$. Fewer than one hand in $180$ requires looking beyond the second relative label.

\section{The big deck, rotation-proof}
\label{sec:bigdeck}
For the $124$-card deck, Berlekamp's recipe is half rotation-proof.
Its first step, which chooses the hidden card by the remainder of the total, is rotation-proof. Place the numbers $0$ to $123$ around a circle and add the same number to every card, going from $123$ back to $0$. The recipe then hides the shifted copy of the old hidden card. (If $t$ cards pass $123$, the total goes up by $t$ modulo $5$, since $124\equiv-1\pmod5$, and so does the position of each card in the sorted hand.)
For example, adding $20$ to every card of the hand $\{7,23,51,88,119\}$ of Section~\ref{sec:ceiling} gives $\{15,27,43,71,108\}$, and the hidden card changes from $88$ to $108$. 
The remaining steps, which encode the hidden card's position, depend on where the numbering starts, so they are not rotation-proof. For example, the shifted hand gets the row $71,27,15,43$, not the shifted row $71,15,43,27$. We replace these steps by a walk around the circle.

\paragraph{A first example.} Before giving the general rule, here is one complete example. Place the numbers $0,1,\ldots,123$ clockwise around a circle and consider the hand $\{4,19,44,78,121\}$.
The sum is $266$, with remainder $1$, so Berlekamp's remainder rule hides the second card in increasing order, namely $19$, and displays $4,44,78,121$. 
The displayed sum is $247$, and $247+3$ is divisible by $5$, so $r=3$. From these four cards alone, the partner can find the $24$ cards that could be hidden. They are in positions $3,8,13,\ldots,118$ when the $120$ unseen cards are listed in increasing order. In terms of labels they are
$$
\begin{aligned}
3,9,14,19,24,29,34,39,45,50,55,60,65,70,75,81,\\
86,91,96,101,106,111,116,122.
\end{aligned}
$$
We call these $24$ cards the \emph{candidates}. Now start at $19$ and walk counterclockwise. The first displayed card met is $4$, and walking clockwise from $4$ back to $19$, the candidates strictly between them are $9$ and $14$. So the count is $2$. 
The assistant chooses $4$ as the lead and shows it first. The order of the other three cards sends the number $2$. Their clockwise order from the lead is $44,78,121$. Listing their six possible orders in dictionary order and numbering them from $0$, order $2$ is $78,44,121$.
The partner sees that the lead is $4$, compares the order shown with the clockwise order $44,78,121$, and recovers the number $2$. She walks clockwise from $4$ with a score of $0$. At the candidate $9$ the score is not yet $2$, so she raises it to $1$, and at $14$ she raises it to $2$. At the next candidate, $19$, the score already equals $2$, so she stops and names $19$.

This example contains the whole idea in its simplest form. In general the lead need not be the first displayed card met, and during the partner's walk other displayed cards may be crossed. The general rule below subtracts $6$ at each such crossing, and the argument after it shows why the same counting idea then works in every case.
\begin{rulebox}{The candidates}
Given the four displayed cards, let $S$ be their sum and let $r\in\{0,1,2,3,4\}$ satisfy $S+r\equiv0\pmod5$. List the $120$ unseen cards in increasing order, with positions numbered $0,\ldots,119$. By Section~\ref{sec:ceiling}, the possible hidden cards occupy positions $r,r+5,\ldots,r+115$. These $24$ cards are the \emph{candidates}.
\end{rulebox}

Every candidate is a possible hidden card. If an unseen card $x$ has $j$ displayed cards below it, its position among the unseen cards is $p=x-j$. If $x$ is a candidate, then $S+x-j=S+p\equiv0\pmod5$, so the total $S+x$ leaves remainder $j$. Now add $x$ to the four displayed cards. In this hand $x$ is in position $j$, so the assistant's rule would hide $x$. Together with Section~\ref{sec:ceiling}, this shows that the candidates are exactly the unseen cards that the assistant's rule would hide if they were added to the four displayed cards. Since Berlekamp's remainder rule is rotation-proof, rotating the displayed cards rotates this set of candidates as well.

\begin{rulebox}{The assistant}
Hide a card using Berlekamp's remainder rule from
Section~\ref{sec:ceiling}. Starting from the hidden card, walk
counterclockwise and meet the four displayed cards one by one. Suppose
a displayed card is the $k$-th one met, and let $c$ be the number of
candidates strictly between that displayed card and the hidden card
when one walks clockwise. Choose the first displayed card for which
$c<6k$. That card is the \emph{lead} and is shown first. Put
$s=c-6(k-1)$. Call the other three displayed cards $A,B,C$ in clockwise
order from the lead. Show them in order number $s$ of
the list
\[
ABC,\quad ACB,\quad BAC,\quad BCA,\quad CAB,\quad CBA,
\]
numbered from $0$.
\end{rulebox}

A lead always exists and $s$ is always one of $0,\ldots,5$. If the first displayed card succeeds, then $k=1$ and $0\leq c<6$. If the first successful card has $k>1$, the preceding displayed card failed and had at least $6(k-1)$ candidates between it and the hidden card. Moving farther counterclockwise cannot decrease that count, so $6(k-1)\leq c<6k$. The fourth displayed card must succeed, because the hidden card is excluded from the strict interval, so $c\leq 23<24$.

\begin{rulebox}{The partner}
The first displayed card is the lead. Call the other three
displayed cards $A,B,C$ in clockwise order from the lead.
Their order in the displayed row determines $s$ from the
same six-order list.
Now walk clockwise from the lead with a score starting at $0$. Whenever a candidate is reached, first check whether the score equals $s$. If it does, stop and name that candidate. If not, add $1$ and continue. Whenever one of the other displayed cards is passed, subtract $6$. The score may become negative.
\end{rulebox}

\begin{figure}[ht]
\centering
\begin{tikzpicture}[scale=1]
  \def\R{2.8}
  \draw[gray!40] (0,0) circle (\R);

  \foreach \cc in {0,...,123}{
    \pgfmathsetmacro\ang{90-\cc*(360/124)}
    \draw[gray!45] (\ang:\R-0.06) -- (\ang:\R+0.06);
  }

  \foreach \cc in {
    0,5,11,16,21,27,32,37,42,47,53,58,
    63,68,73,78,83,93,98,103,108,113,118
  }{
    \pgfmathsetmacro\ang{90-\cc*(360/124)}
    \fill[blue!60!black] (\ang:\R) circle (1.9pt);
  }

  \pgfmathsetmacro\ang{90-88*(360/124)}
  \node[
    star,star points=5,fill=red!75!black,inner sep=1.8pt
  ] at (\ang:\R) {};
  \node[font=\small\bfseries,red!75!black]
    at (\ang:\R+0.5) {88};

  \foreach \cc in {7,23,51,119}{
    \pgfmathsetmacro\ang{90-\cc*(360/124)}
    \fill[black] (\ang:\R) circle (3.3pt);
    \node[font=\small\bfseries] at (\ang:\R+0.5) {\cc};
  }

  \pgfmathsetmacro\sa{90-7*(360/124)}
  \pgfmathsetmacro\ea{90-88*(360/124)}
  \draw[
    -{Stealth[length=3mm]},
    line width=1.2pt,green!45!black
  ]
    (\sa:\R-0.45)
    arc[start angle=\sa,end angle=\ea,radius=\R-0.45];

  \node[font=\small,green!35!black,align=center]
    at (0.3,0.35)
    {partner walks\\clockwise from the\\lead card $7$};

  \node[font=\scriptsize,align=left,anchor=west]
    at (\R+0.9,-1.9) {%
      \tikz\fill[black] (0,0) circle (3pt); displayed\\
      \tikz\fill[blue!60!black] (0,0) circle (2pt); candidate\\
      \tikz\node[
        star,star points=5,fill=red!75!black,inner sep=1.5pt
      ]{}; hidden
    };
\end{tikzpicture}
\caption{The hand $\{7,23,51,88,119\}$ on the $124$-card circle.
Berlekamp's remainder rule hides $88$. The candidates are every fifth undisplayed card.
Walking counterclockwise from $88$, the $51$ has $7$ candidates strictly between it and $88$ (not fewer than $6$), the $23$ has $12$ (not fewer than $12$), and the $7$ has $15$ (fewer than $18$). So $7$ leads, $s=15-12=3$, and the row is $7,51,119,23$. On the walk clockwise from $7$, the first candidate reached with score $3$ is $88$.
}
\label{fig:circle}
\end{figure}
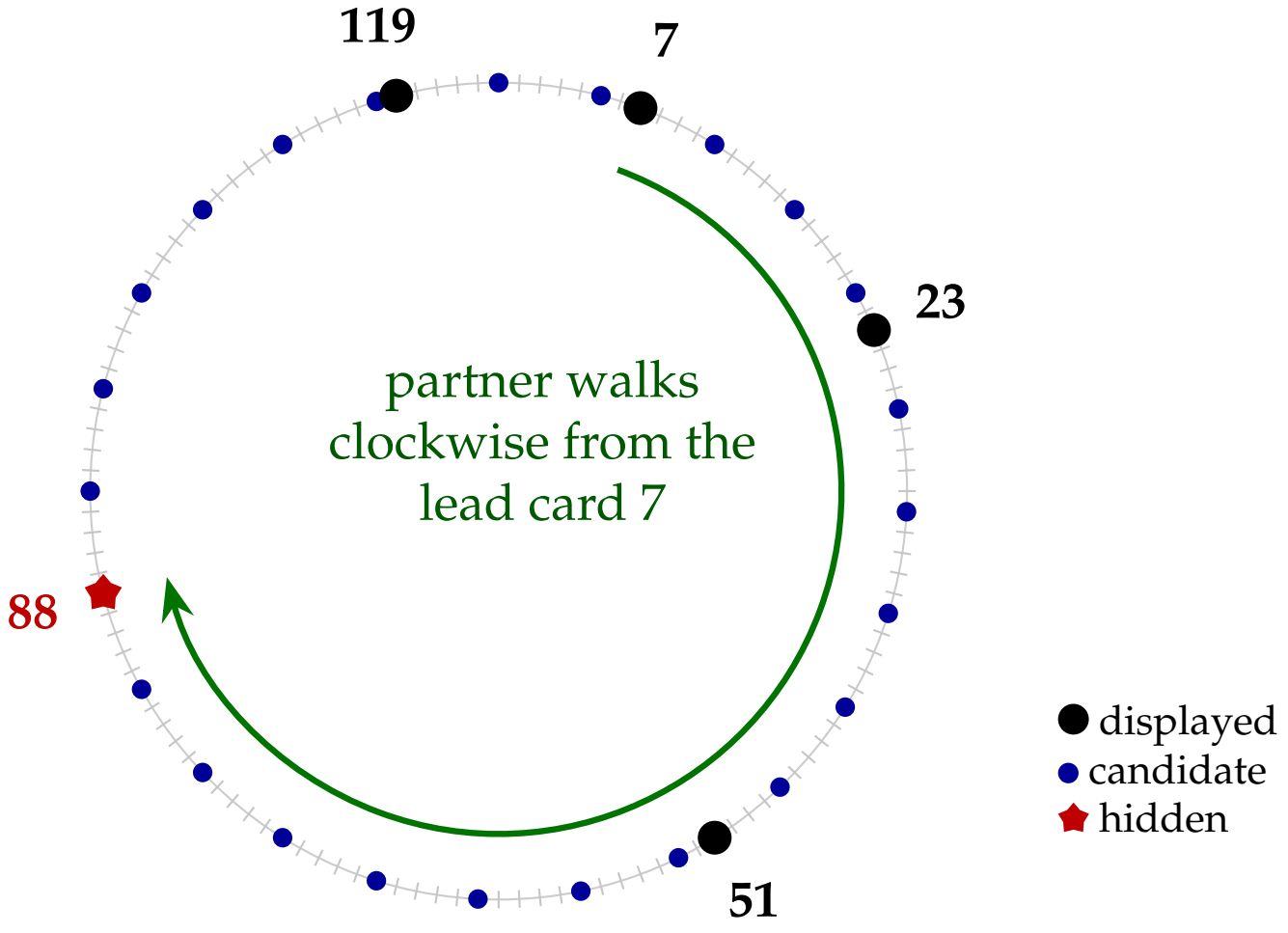

For the example in Figure~\ref{fig:circle}, walking counterclockwise
from $88$ meets $51$ ($k=1$, $c=7\not<1\cdot6$), then $23$ ($k=2$,
$c=12\not<2\cdot6$), then $7$ ($k=3$, $c=15<3\cdot6$). Thus $7$ is the lead and
$s=15-12=3$. 
In clockwise order from the lead, the other three cards are $23,51,119$, and order number $3$ of the list is $51,119,23$, so the row is $7,51,119,23$.
Starting from $7$, the partner counts candidates
$11,16,21$, bringing the score to $3$. She does not stop at $21$,
because the score was $2$ before counting it. Passing $23$ reduces it
to $-3$. Counting $27,\ldots,47$ brings it to $2$, and passing $51$ to
$-4$. Counting $53,\ldots,83$ brings it to $3$. At the next candidate,
$88$, the score already equals $3$, so she stops and names $88$. 

The score rule reverses the assistant's construction. If the lead was the $k$-th displayed card met counterclockwise from the hidden card, then $k-1$ other displayed cards lie between the lead and the hidden card clockwise. Immediately before counting the hidden candidate, the score has increased by $c$ and has been reduced by $6$ exactly $k-1$ times, so it equals $c-6(k-1)=s$ and the partner stops.

\paragraph{Why the partner does not stop too soon.}
Consider an intervening displayed card that was the $i$-th card met counterclockwise from the hidden card, with $c_i$ candidates strictly between it and the hidden card. It failed the assistant's test, so $c_i\geq6i$. For a candidate in the stretch just before this displayed card on the clockwise walk, at most $c-c_i-1$ candidates have already been counted and exactly $k-1-i$ displayed cards have been passed, so the score is at most
$$
c-c_i-1-6(k-1-i)=s+6i-c_i-1\leq s-1.
$$
In the last stretch, after all $k-1$ of these displayed cards, the score at any candidate before the hidden card is at most $c-1-6(k-1)=s-1$. This also covers $k=1$. So the first stop is the hidden card.

Now return to the hand $\{15,27,43,71,108\}$ from the beginning of this section, the hand of Figure~\ref{fig:circle} shifted by $20$. The new method hides $108$ and displays $27,71,15,43$, which is the row $7,51,119,23$ with $20$ added to each card, whereas Berlekamp's encoding gave $71,27,15,43$.
From the lead $27$, the other three cards in clockwise order are $43,71,15$. The row shows them as $71,15,43$, which is again order number $3$, so $s=3$.

Once the candidates are known, the remaining steps use only the order of the cards around the circle and the number of candidates between cards, and a rotation changes neither. Therefore the whole method is rotation-proof, and the performers may start the numbering anywhere, as long as they both start at the same place.
As in Section~\ref{sec:nostart}, no rotation other than the identity carries a hand onto itself, because its cycle length would have to divide both $124$ and $5$. So the hands fall into $\binom{124}{5}/124=1{,}815{,}726$ rotation families of $124$ hands each, and checking one hand from each family is enough.
A computer check of all $225{,}150{,}024$ hands confirms that the partner always names the hidden card and that shifting a hand shifts its row. The same walk works for the smaller tricks of Section~\ref{sec:variations}, with $5$ replaced by the hand size, $3$ or $4$, and $6$ by the number of ways to order the cards after the lead, $1$ or $2$. 
These versions are rotation-proof too, because $8\equiv-1\pmod3$ and $27\equiv-1\pmod4$ play the role of $124\equiv-1\pmod5$, and they were checked on all $56$ and $17{,}550$ hands.

A related circular counting device appears in a method posted by \emph{ThirdParty} on the xkcd puzzle forum and described by Craig Gidney \cite{Gidney}.

\section{Ideas for performers}

\begin{enumerate}
\renewcommand{\labelenumi}{(\roman{enumi})}

\item \textbf{Use your own suit cycle.}
The rotation-proof method needs a cyclic order of the suits,
not a linear one, and any cycle works as long as both performers
use the same order and direction. A performer whose stack uses
a repeating suit cycle can use that same cycle.

\item \textbf{Let the spectator set the clock.}
Because the method has no starting point, nothing depends on
which rank is called ``first.'' The spectator may choose which
rank goes at the top of a drawn clock, while the cyclic order
and clockwise direction stay unchanged. This gives the spectator
another choice without changing the method.

\item \textbf{The partner does not break ties between pairs.}
Both performers use the relative ranking of the signaling cards
on every hand. The partner reads their order and counts forward
from the base as before. Choosing the pair is the assistant's
task. Only about $15\%$ of hands require further tie-breaking
after the smallest pair distance has been found.

\item \textbf{The pocket version}
with eight cards gives a small setting in which to teach and
practice the method. It works with any eight distinct cards:
agree on an order, assign them the values $1,\ldots,8$, and
use these assigned values in the arithmetic.

\item \textbf{The heckler challenge}
makes a closing talking point: with a full deck, no method using
only the displayed row can guarantee success for every hand if
a spectator may swap two neighboring cards or leave the row
unchanged, without telling the partner. With $29$ cards, however,
the challenge can be met.
\end{enumerate}

\section*{Acknowledgments}

The authors used Claude Opus 5.5 (Anthropic) to assist with the literature search, with writing the computer programs that check the methods described here, with writing the TikZ code for the figures, and with the preparation of the text. The authors take full responsibility for the content of the article.

\end{document}